\documentclass{article}
\usepackage{amsmath,amssymb,enumerate,bbm,theorem}

\newcommand{\email}[1]{{\textit{Email:} \texttt{#1}}}
\newcommand{\nequiv}{\not\equiv}
\newcommand{\nin}{\not\in}
\newcommand{\tmem}[1]{{\em #1\/}}
\newcommand{\tmop}[1]{\ensuremath{\operatorname{#1}}}
\newcommand{\tmsamp}[1]{\textsf{#1}}
\newcommand{\tmstrong}[1]{\textbf{#1}}
\newcommand{\tmtexttt}[1]{{\ttfamily{#1}}}
\newenvironment{enumeratenumeric}{\begin{enumerate}[1.] }{\end{enumerate}}
\newenvironment{proof}{\noindent\textbf{Proof\ }}{\hspace*{\fill}$\Box$\medskip}
\newtheorem{definition}{Definition}
{\theorembodyfont{\rmfamily}\newtheorem{example}{Example}}
\newtheorem{lemma}{Lemma}
\newtheorem{proposition}{Proposition}
{\theorembodyfont{\rmfamily}\newtheorem{remark}{Remark}}
\newtheorem{theorem}{Theorem}

\begin{document}

\title{Trees, valuations and the Green-Lazarsfeld set.}\author{Thomas
Delzant\thanks{\email{delzant@math.u-strasbg.fr}}\\
D\'epartement de Math\'ematiques, Universit\'e de Strabourg\\
7 rue R. Descartes F-67084 Strasbourg}\maketitle

\tmtexttt{}\section{Introduction.}

The aim of this paper is the study of the relationship between two objects,
the Green-Lazarsfeld set and the Bieri Neumann Strebel invariant, which appear
simultaneously in 1987 ($\left[ \tmop{GL} \right], \left[ \tmop{BNS}
\right])$. Let us recall some basic definitions.

Let $\Gamma$ be a finitely generated group, and $K$ be a field. A
$1$-character $\chi$ is an homomorphism from $\Gamma$ to $K^{\ast}$ ; in this
article we will only consider 1-characters, and call them characters. A
character $\chi$ is called exceptional if $H^1 (\Gamma, \chi) \neq 0$, or more
geometrically if $\chi$ can be realized as the linear part of a \ fixed point
free affine action of $\Gamma$ on a $K$-line.

The set of exceptional characters, $E^1 (\Gamma, K)$ is a subset of the
abelian group $\tmop{Hom} (\Gamma, K^{\ast})$, and our aim is to understand
its geometry, in particular if $\Gamma$ is the fundamental group of a K\"ahler
manifold.

Motivated by the pioneering work of M. Green and R. Lazarsfeld, algebraic
geometers studied the case where $K =\mathbbm{C}$ is the field of complex
numbers, and $\Gamma = \pi_1 (X)$ is the fundamental group of a projective or
more generally a K\"ahler manifold. In this case, the geometry of $\tmop{Hom}
(\Gamma, K^{\ast})$is well understood : it is the union of a finite set, made
up with torsion characters, and a finite set of translates of subtori. This
result has been conjectured by A. Beauville and F. Catanese, ($\left[
\right.$Be$\left. \right])$ proved by C. Simpson $\left[ \tmop{Si} 2 \right]$)
for projective manifolds, and extended by F. Campana to the K\"ahler case (see
$\left[ \tmop{Ca} \right]$ for a detailed introduction). The main tools used
by C. Simpson were the flat hyper-K\"ahler structure of $\tmop{Hom} (\Gamma,
\mathbbm{C}^{\ast})$ and the \ Schneider-Lang theorem in transcendence theory.
Another proof, model theoretic, has been proposed by R. Pink and D. Roessler
$\left[ \tmop{PR} \right]$.

The definition of a exceptional class in the sense of Bieri Neumann Strebel is
easier to explain in the case of an integral cohomology class (an element of
$H^1 (\Gamma, \mathbbm{Z})$). Such a class is {\tmem{exceptional}} if it can
be realized as the translation class of a parabolic, non loxodromic action of
$\Gamma$ in some tree.

The link between these two notions, explained in the next paragraph, can be
sketched as follows. Let $\chi$ be a Green-Lazarsfeld character of $\Gamma$.
Suppose $\chi (\Gamma)$ is not contained in the ring of algebraic integers of
$K$. There exists a discrete non archimedian valuation on $K$ such that $v
\circ \chi$ is a non trivial homomorphism to $\mathbbm{Z}$. It appears that $v
\circ \chi$ is an exceptional class in the sense of Bieri Neumann Strebel.
More precisely, on can find a parabolic action of $\Gamma$ on the Bruhat Tits
tree of $K_v$, \ with translation length $v \circ \chi .$

Due to the work of C. Simpson $\left[ \tmop{Si} 3 \right],$ M. Gromov and R.
Shoen $\left[ \tmop{GS} \right]$, exceptional cohomology classes on K\"ahler
manifold are well understood (see also $\left[ \tmop{De} \right]$ for a
detailed study of the BNS invariant of a K\"ahler group). Let $X$ be a
K\"ahler manifold, and $\omega$ \ an exceptional class ; there exists a
holomorphic map $F$ from $X$ to a \ hyperbolic Riemann orbifold $\Sigma$such
that $\omega$ belongs to $F^{\ast} H^1 (\Sigma, \mathbbm{Z})$. Recall that a
complex 2-orbifold $\Sigma$ is a Riemann surface $S$ marked by a finite set of
marked points $\left\{ (q_1, m_1) \ldots, (q_n, m_n) \right\} $, where the
$m'_i s$ are integers $\geqslant 2$. A map $F : X \rightarrow \Sigma$ is
called holomorphic if it is holomorphic in the usual sense, and for every
$q_i$ the multiplicity of the fiber $F^{- 1} (q_i)$ is divisible by $m_i .$
The main result of this paper is a description of the (generalized) Green
Lazarsfeld set of $\pi_1 (X)$ in terms of the finite list of its fibrations on
hyperbolic 2-orbifolds.

\paragraph{Theorem}{\tmem{Let $\Gamma$ be the fundamental group of a K\"ahler
manifold X, $(F_i, \Sigma_i)_{1 \leqslant i \leqslant n}$ the family of
fibration of $X$ over hyperbolic $2$-orbifolds. Let $K$ be a field of
characteristic $p$ (if $p = 0$, $K =\mathbbm{C})$, $\bar{F}_p \subset K$ the
algebraic closure of $F_p$ in $K$. Then $E^1 (\Gamma, K)$ is the union of a
finite set of torsion characters (contained in $E^1 (\Gamma, \bar{F}_p$) if $p
> 0$) and the union}} \ $\bigcup_{1 \leqslant i \leqslant n} F_i^{\ast} E^1
(\pi^{\tmop{orb}}_1 (\Sigma_i), K^{\ast}) .$

\paragraph{Remarks}a) Let $\Sigma = (S ; (q_i, m_i)_{1 \leqslant i \leqslant
n})$ a hyperbolic \ 2-orbifold, and \ $\Gamma = \pi_1^{\tmop{orb}} (\Sigma)$
its fundamental group. Then, by a simple computation (see prop.22), \ one
checks that $E^1 (\pi^{\tmop{orb}}_1 (\Sigma_{}), K^{\ast}) = \tmop{Hom}
(\pi_1^{\tmop{orb}} (\Sigma, K^{\ast})) = (K^{\ast})^{2 g} \times \Phi$, where
$\Phi$ is a finite abelian group, unless $g = 1$ and for all $i$, $m_i \nequiv
0 (\tmop{char} K) .$ If $g = 1$ and for all $i$ $m_i \nequiv 0 (\tmop{char}
K)$, \ $E^1 (\pi^{\tmop{orb}}_1 (\Sigma_{}), K^{\ast})$ is finite, made of
torsion characters.

In every cases, the Green Lazarsfeld set is the union of a finite set of
torsion characters and a finite set of abelian groups which are translates of
tori ; this is our generalization of Simpson's theorem.

b) The main tool used by Simpson to prove his theorem $\left[ \tmop{Si} 2
\right] $was the study of algebraic triple tori ; if $\tmop{char} K \neq 0$ no
such a structure is available. Our proof furnishes a geometric (i.e. non
arithmetic) alternative to Simpson's proof in the case of characteristic 0. In
fact, in this case ($\tmop{car} K = 0$) our method proves that \ $E^1 (\Gamma,
K)$ is made with a finite set of {\tmem{integral}} characters (in the sense of
Bass $\left[ \right.$Ba$\left. \right]$), and the union $\bigcup_{1 \leqslant
i \leqslant n} F_i^{\ast} \tmop{Hom} (\pi^{\tmop{orb}}_1 (\Sigma_i),
K^{\ast})$ \ ; the conclusion follows from the study of the the absolute value
$| \chi |$ of exceptional characters, which was already done by \ A.
Beauville$ \left[ \tmop{Be} \right]$.

\paragraph{}

In a recent preprint, \ $\left[ \tmop{CS} \right],$C. Simpson and K. Corlette
study the variety of characters of a K\"ahler group $\Gamma$, \
$\tmop{Hom}^{\tmop{ss}} (\Gamma, \tmop{PSl} (2, \mathbbm{C}) / \tmop{PSL} (2,
\mathbbm{C})$ from a very similar point of view ; they prove in particular
that a Zariski dense representation of a K\"ahler group which is not
{\tmem{integral}} in the sense of Bass factorizes through a fibration over a
hyperbolic 2-orbifold. Their proof is based on the same idea as ours : if a
representation $\rho$ is not integral, there exists a valuation on the field
generated by $\rho (\Gamma)$ such that the action of $\Gamma$ on the Bruhat
Tits building is non elementary. The conclusion follows by applying the theory
of Gromov Shoen on harmonic maps with value in a tree. Using Simpson's work on
Higgs bundles they prove further a rigid representation come from a complex
variation of Hodge structure.

In paragraph 2, we explain the relationship between the Green-Lazarsfeld and
Bieri-Neumann-Strebel invariants ; in paragraph 3 we study the
Green-Lazarsfeld set of a metabelian group : a finiteness result on this set
is established. These two paragraphs are purely group theoretic, and no
K\"ahler structure is mentioned. In the paragraph 4 we prove the main result.

{\tmstrong{Acknowledgments.}} I would like to thank R. Bieri for very helpful
discussions on the structure of metabelian groups, and for explaining me his
paper $\left[ \tmop{BG} \right]$ with J. Groves, and F. Campana for his
interest and comments.

\section{From an \ affine action on a line to a parabolic action on a tree. }

\subsection{Affine action on the line : the Green-Lazarsfeld set}

Let $K$ be a field. The affine group of transformation of a $K$-line,
$\tmop{Aff}_1 (K),$ is isomorphic to $K^{\ast} \ltimes K$. We identify this
group with the set of upper triangular $(2, 2)$matrices
$\left(\begin{array}{cc}
  \ast & \ast\\
  0 & 1
\end{array}\right)$ with values in $K$.

Let $\Gamma$ be a finitely generated group. An affine action of $\Gamma$ on
the line is a morphism $\rho : \Gamma \rightarrow \tmop{Aff}_1 (K)$. One can
write $\rho (g) = \left(\begin{array}{cc}
  \chi (g) & \theta (g)\\
  0 & 1
\end{array}\right) .$ The linear part of $\rho$ is an homomorphism $\chi :
\Gamma \rightarrow K^{\ast}$. Its translation part $\theta : \Gamma
\rightarrow K$ is a 1-cocycle of $\Gamma$ with value in $\chi$, i.e. a
function which satisfies $\theta (g h) = \theta (g) + \chi (g) \theta (h)$.
The representation $\rho$ is conjugate to a diagonal representation if and
only if $\rho (\Gamma)$ fixes a point $\mu \in K$, or equivalently if and only
if there exists a $\mu \in K$ such that $\theta (g) = \mu (- 1 + \chi (g))$ is
a coboundary.

\begin{definition}
  A character $\chi \in \tmop{Hom} (\Gamma, K^{\ast})$is exceptional if it can
  be realized as the linear part of a fixed point free affine action of
  $\Gamma$ on the line, i.e if $H^1 (\Gamma, \chi) \neq 0$. The set of
  exceptional characters $E^1 (\Gamma, K)$ is called the Green-Lazarsfeld set
  of $\Gamma .$
\end{definition}

\subsection{Parabolic action on a tree : the Bieri Neumann Strebel invariant.}

Let \ $T$ be a simplicial tree. We endow $T$ with its natural simplicial
metric, and think of $T$ as a complete geodesic space. Let us recall the
definitions of the boundary of $T$, and of the Busemann cocyle associated to a
point in this boundary.

A ray in $T$ is an isometric map $r : \left[ a, + \infty \left[ \rightarrow T
\right. \right.$. Two rays $r : \left[ a, + \infty \left[ \rightarrow T
\right. \right.$, $s : \left[ b, + \infty \left[ \rightarrow T \right.
\right.$ are equivalent (or asymptotic) if they coincide after a certain time
: there exists $a', b'$ s.t. for all $t \geqslant 0$ $r (a' + t) = s (b' +
t)$. The boundary of $T$, denoted $\partial T,$ is the set of equivalence
classes of rays. If $\alpha \in \partial T$ and $r : \left[ a, + \infty \left[
\rightarrow T \right. \right.$ represents $\alpha$, for every point $x$, the
function $t \rightarrow d (x, r (t)) - t$ is eventually constant. Its limit
$b_r (x)$ is called the Busemann function of $r$. If $s$ is equivalent to $r$,
the difference $b_r - b_s$ is a constant.

\begin{definition}
  (Busemman cocyle) . Let $\Gamma$ be a group acting on $T$, and $\alpha \in
  \partial T$. If $\Gamma$ fixes $\alpha$, \ one \ define an homomorphism, the
  \ Busemann cocyle, by the formula :

  \ \begin{center}
    $\omega_{\alpha} : \Gamma \rightarrow \mathbbm{Z}$
  \end{center}
  
  $\omega_{\alpha} (g) = b_r \circ g - b_r$
\end{definition}

\begin{definition}
  (Exceptional classes) The action of $\Gamma$ is called parabolic if it fixes
  some point at infinity. It is called exceptional if fixes a unique point at
  infinity, and if the associated Busemann cocycle is not trivial. A class
  $\omega \in H^1 (\Gamma, \mathbbm{Z})$ is exceptional if it can be realized
  as the Busemann cocycle of an exceptional action of $\Gamma$ in some tree.
  The set of exceptional classes is denoted $\mathcal{E}^1 (\Gamma,
  \mathbbm{Z}) .$
\end{definition}

\begin{remark}
  A {\tmem{topological }}definition of a exceptional class can also be given,
  in the case where $\Gamma$ is finitely presented. \ Let $\Gamma = \pi_1
  (X),$where $X$ is a compact manifold, and let $\omega$ be some class in $H^1
  (\Gamma, \mathbbm{Z}) .$ One represents $\omega$ by a closed 1-form $w$ on
  $X$ and consider a primitive $F : \tilde{X} \rightarrow \mathbbm{R}$ \ of
  the lift of $w \tmop{to} \tmop{the} \tmop{universal}$ cover of $X.$ Then
  $\omega$ is exceptional iff $F \geqslant 0$ has several components on which
  $F$ is unbounded (see $\left[ \tmop{Bi} \right], \left[ \tmop{Le} \right],
  \left[ \tmop{Bro} \right])$.
\end{remark}

\begin{remark}
  The notion of an exceptional class, defined by Bieri Neumann Strebel and
  studied by several authors, in particular $\left[ \tmop{Bro} \right], \left[
  \tmop{Le} \right],$is more general : it concerns homomorphism with value in
  $\mathbbm{R}$ and can be defined along the same lines, using
  $\mathbbm{R}$-trees instead of combinatorial trees. Our point of view is
  that of Brown ; it is interesting to remark that $\left[ \tmop{Bro} \right],
  \left[ \tmop{BNS} \right] $and $\left[ \tmop{GL} \right]$ are published in
  the same issue of the same journal, but apparently nobody remarked that
  $\left[ \tmop{Bro} \right]$ and $\left[ \tmop{GL} \right]$ studied the same
  object from a different point of view. This remark justify the choice of our
  title. 
\end{remark}

\subsection{Discrete valuations and \ Bruhat-Tits trees.}

In this paragraph we fix a field $K.$ Let $v : K^{\ast} \rightarrow
\mathbbm{Z}$ be a discrete non archimedian valuation on $K$. Bruhat and Tits
$\left[ \tmop{BT} \right] $constructed a tree $T_v$ with an action of
$\tmop{PGL} (2, K)$. One should think of the action of $\tmop{PGL} (2, K)$ of
$T_v$ as an analogue of the action of $\tmop{PGL} (2, \mathbbm{C})$ on the
hyperbolic space of dimension 3 ; \ we recall below some basic facts about
this action (see $\left[ \tmop{Se} \right]$ for a detailed study).

Let $O_v \subset K$ denote the valuation ring $v \geqslant 0$. The
{\tmem{vertices}} of $T_v$ are the homothety classes of $O_v -$lattices, i.e.
free $O_v -$modules of rank 2, in $K^2$ . The {\tmem{boundary}} of this tree
is the projective line $P^1 (^{} \bar{K}_v)$ over the $v$-completion of $K.$

By the general theory of lattices, if $\Lambda, \Lambda'$ are two lattices,
one can find a $O_v$-base of $\Lambda$ such that, in this base, $\Lambda'$ is
generated by $(t^a, 0)$ and $(0, t^b)$ for some $t$ with $v (t) = 1$ ; hence
up to homothety by $(1, 0)$ and $(0, t^n)$, for $n = b - a.$ Then the distance
between $\Lambda,$ and $\Lambda'$ is $| n |$, and the segment between
$\Lambda$ and $\Lambda'$ is the set of lattices generated by $(1, 0)
\tmop{and} (0, t^k), k = 1, n$. More generally if $l, l'$ are two different
lines in $K^2$, considered as points in $\partial T_v$, the geodesic from $l$
to $l'$ is the set of product of lattices in $l$ and $l' .$

The matrix $g_u = \left(\begin{array}{cc}
  1 & u\\
  0 & 1
\end{array}\right)$ fixes the lattice $\Lambda_n $ generated by $(1, 0)$ and
$(0, t^n)$ for $n \leqslant v (u) .$ The matrix $g_u = \left(\begin{array}{cc}
  t^n & u\\
  0 & 1
\end{array}\right)$ transforms $\Lambda_m \tmop{to} \Lambda_{m + n}$ if $m + n
\leqslant v (u) .$

Acting on $T_v $ the Borel sub-group $\left(\begin{array}{cc}
  \ast & \ast\\
  0 & 1
\end{array}\right)$ is {\tmem{parabolic}} : it fixes an end of $T_v$ (namely
the line generated by the first basis vector), but neither a point of $T_v$
nor a pair of points of $\partial T_v .$

The Busemann cocyle of this parabolic subgroup is $b \left( \right. \left.
\left(\begin{array}{cc}
  \alpha & \beta\\
  0 & 1
\end{array}\right) \right) = v (\alpha) .$

The relation between the Green-Lazarsfeld set and the Bieri-Neumann-Strebel
invariant is now simple to explain. \

\begin{proposition}
  Let $\chi \in H^1 (\Gamma, K^{\ast})$. Suppose that $\chi \in E^1 (\Gamma,
  K^{\ast})$ and let $\theta \in H^1 (\Gamma, \chi) \neq 0$. Let $\rho :
  \Gamma \rightarrow \tmop{Gl} (2, K)$ be defined by $\rho (g) =
  \left(\begin{array}{cc}
    \chi (g) & \theta (g)\\
    0 & 1
  \end{array}\right) .$ If $v \circ \chi \in H^1 (\Gamma, \mathbbm{Z})$ is not
  0, $\rho$ is an exceptional action on $T_v$.
\end{proposition}

{\tmstrong{Proof.}} \ By construction the action of $\Gamma$ on $T_v$ fixes a
point at infinity. It contains an hyperbolic element as $v \circ \chi \neq 0$,
but the action cannot fix a line : the other point in the boundary $P^1
(^{^{^{}}} \bar{K}_v)$ would be fixed by the group $\Gamma$, and $\rho$ would
be conjugate to diagonalizable action. The orbit of any point of $\Gamma$ is
therefore a minimal tree which is not a line. $\Box$

\section{Metabelian groups}

If $\Gamma$ is a group, let $\Gamma' = \left[ \Gamma, \Gamma \right]$ its
derived group. Recall that a group is{\tmem{ metabelian }}if $\Gamma'$ is
abelian, or $\Gamma^2 = (\Gamma')'$ is trivial. If $\Gamma$ is a f.g. group,
$\Gamma / \Gamma^2$ is metabelian.

\subsection{The Green-Lazarsfeld set of a metabelian group.}

If $K$ is a field, the Green-Lazarsfeld set $E^1 (\Gamma, K)$ of the group
$\Gamma$ only depends on its metabelianized $\Gamma / \Gamma^2$ \ as it only
depends of the set of representation of $\Gamma$ in the metabelian group
$\tmop{Aff}_1 (K) = K^{\ast} \ltimes K$ .

Let $\Gamma$ be a metabelian group. We write $1 \rightarrow \left[ \Gamma,
\Gamma \right] \rightarrow \Gamma \rightarrow Q \rightarrow 1$, where $Q =
\Gamma / \left[ \Gamma, \Gamma \right]$ is the abelianized group, and $\left[
\Gamma, \Gamma \right]$ is abelian. As an abelian group, $M = \left[ \Gamma,
\Gamma \right]$ is not necessary f.g, however we can let $Q$ acts on $\left[
\Gamma, \Gamma \right]$ by conjugation, so that $M$ can be promoted as a
$\mathbbm{Z}Q$ module. The following fact is basic and well-known.

\begin{lemma}
  The module $M$ is finitely generated as a $\mathbbm{Z}Q$ module.
\end{lemma}

If $g_1, \ldots \ldots g_r$ are generators of $\Gamma$, the commutators $h_{i
j} = \left[ g_i, g_j \right]$ generate $\text{$\left[ \Gamma, \Gamma
\right]$}$ as a $\mathbbm{Z}Q$ module : if $\left[ g, h \right]$ if $h = a b$
\ we have $\left[ g, h \right] = \left[ g, a b \right] = g a g^{- 1} a^{- 1} a
g b g^{- 1} b^{- 1} a^{- 1} = \left[ g, a \right] a \left[ g, b \right] a^{-
1} = \left[ g, a \right] a_{\ast} \left[ g, b \right]$, and the result follows
by induction. $\Box$

\begin{theorem}
  Let $\Gamma$ be a finitely generated group. Given a prime number p (p might
  be 0), there exists a finite number of fields $K_{\nu}$ of characteristic
  $p$ and of finite transcendence degree over $F_p $ (if $p = 0$, set $F_p
  =\mathbbm{Q})$ and characters $\xi_{\nu} : \Gamma \rightarrow
  K_{\nu}^{\ast}$ such that :
  \begin{enumeratenumeric}
    \item $H^1 (\Gamma, \xi_{\nu}) \neq 0$, i.e. $\xi_{\nu} \in E^1 (\Gamma,
    K_{\nu})$
    
    \item If $K$ is a field of characteristic $p$ and $\chi \in E^1 (\Gamma,
    K)$ a Green-Lazarsfeld character, then there exists an index $\nu$ s.t.
    $\ker \chi \supset \ker \xi_{\nu} .$
  \end{enumeratenumeric}
\end{theorem}

\begin{proof}
  Let $F_p$ be the field with $p$ elements \ and \ $F_p \left[ Q \right]$ the
  group ring of $Q$ with $F_p$ coefficients. Let $M_p = \left[ \Gamma, \Gamma
  \right] \otimes F_p, \mathcal{J} \subset F_p \left[ Q \right]$ the
  annihilator of $M_p$, and $A = F_p \left[ Q \right] /\mathcal{J}$. As $Q$ is
  a finitely generated abelian group, isomorphic to $\mathbbm{Z}^r \times
  \Phi$, with $\Phi$ finite abelian, \ $A$ is a noetherian ring. Thus $A$
  admits a finite number of {\tmem{minimal}} prime ideals
  $(\mathfrak{p}_{\nu})_{1 \leqslant \nu \leqslant \nu_0}$. Let $k_i$ be the
  field of fraction of $A /\mathfrak{p}_i$, and $\xi_i$ be the natural
  character $\Gamma \rightarrow Q \rightarrow A /\mathfrak{p}_i \rightarrow
  k_i .$ Up to re-ordering the list of these ideals, we may assume that for $1
  \leqslant i \leqslant \nu_1, H^1 (\Gamma, \xi_i) \neq 0$.

  The theorem 8 is a consequence of the following :
  
  \begin{lemma}
    Let $\chi \in E^1 (\Gamma, K)$ be an exceptional character, $\chi \neq 1$,
    and let $\mathfrak{p}$ be a minimal prime ideal contained in $\ker \chi$.
    Then, the character $\xi_{\mathfrak{p}}$ belongs to $E^1 (\Gamma,
    k_{\mathfrak{p}}^{\ast})$, i.e. $H^1 (\Gamma, \xi) \neq 0$.
  \end{lemma}

  Let$M_{\mathfrak{p}} = M \otimes A_{\mathfrak{p}}$, and \ $M_0 = M \otimes_A
  K = M_{\mathfrak{p}} /\mathfrak{p}M_p $. Note that $M_0$ is a finitely
  generated $k_{\mathfrak{p}}$ vector space, on which $\Gamma$ acts \ by
  homotheties : the action of $g$ is the homothety of ratio $\xi (g)$. Let
  $\pi : \left[ \Gamma, \Gamma \right] \rightarrow M_0$ the canonical map. We
  shall prove that $H^1 (\Gamma, M) \neq 0$.

  For some $g_0 \in \Gamma, \xi (g_0)$ is not 1 (as an element of
  $k_{\mathfrak{p}}$) : if not $\Gamma = \ker \xi_p$ so $\chi = 1$.

  The map $\Gamma \rightarrow M_0$ defined by $c (g) = \pi (g_0 g g^{- 1}_0
  g^{- 1})$ satisfies $c (g h) = \pi (g_0 g h g^{- 1}_0 h^{- 1} g^{- 1}) = \pi
  (g_0 g g_0^{- 1} g^{- 1}) + \pi (g g_0 h g^{- 1}_0 h^{- 1} g^{- 1}) = c (g)
  + \xi (g) \pi (g_0 h g^{- 1}_0 h^{- 1}) = c (g) + \xi (g) c (h)$. Therefore
  $c$ is a 1-cocycle of $\Gamma$ with value in $M$.

  Let us prove, by contradiction, \ that the cohomology class of $c$ is not
  $0$. \

  For every $m \in M_0$, $c (m) = (\xi (g_0) m - m) = (\xi (g_0) - 1) m$. \
  If $c = 0$, as $\xi (g_0) \neq 1,$ then $M_0 = 0$. But if \ $M_0 = 0,
  M_{\mathfrak{p}} /\mathfrak{p}M_p = 0,$i.e. $\mathfrak{p}M_{\mathfrak{p}} =
  M_{\mathfrak{p}}$, and $M_{\mathfrak{p}} = 0$ by the Nakayama lemma
  ($\mathfrak{p}$ is the unique maximal ideal \ of{\tmsamp{
  $A_{\mathfrak{p}})$}}, i.e. $M =\mathfrak{p}M$. But $\mathfrak{p} \subset
  \ker \chi$, so \ this would implies that \ $M \otimes_A K = 0$ and $H^1
  (\Gamma, \chi_{}) = 0$.

  If this cocyle is a coboundary we could find some $m \in M_0$ s.t. $c (g) =
  (1 - \xi (g)) m,$ but $c (g_0) = 0,$ and $\xi (g_0) \neq 1$, so $c$ would be
  $0$.

  In order to prove lemma 9, we see that, for every linear map $l = M_0
  \rightarrow K$, $l \circ c$ is a non trivial $1 - \tmop{cocycle} .$

  This proves theorem 8.

\end{proof}

\begin{remark}
  The previous proof is a combination of arguments by $\left[ \tmop{BG}
  \right]$ and $\left[ \tmop{Bre} \right]$. In their remarkable paper R. Bieri
  and J. Groves describe the BNS invariant of a metabelian group in terms of
  the finite set of field $k_{\nu} $ and characters $\xi_{\nu}$ for a finite
  set of primes $p$ (the primes $p$ for which $\left[ \Gamma, \Gamma \right]$
  has $p$-torsion). For every such a field and every valuation on it, $v \circ
  \xi_{\nu}$ is exceptional. This provide a map from the cone of valuations on
  $k_{\nu} $to the BNS set. This set turns out to be the union of the images
  of these cones. In $\left[ \tmop{Bre} \right],$ Breuillard proves along the
  same lines, that a metabelian not virtually nilpotent group admits a non
  trivial affine action.
\end{remark}

\section{Fundamental groups of K\"ahler manifolds.}

\subsection{Fibering a K\"ahler manifold.}

\ \ \ \

For the general study of orbifolds and their fundamental groups, we refer to
W. Thurston $\left[ \tmop{Th} \right]$ chap. 13. Complex $2$-orbifolds are
2-orbifolds with singularities modeled on the quotient of the unit disk by the
action of $\mathbbm{Z}/ n\mathbbm{Z}.$ The usefulness of this notion in our
context of (fibering complex manifolds to Riemann surfaces) has been pointed
out by C. Simpson $\left[ \tmop{Si} 1 \right] .$

\begin{definition}
  Complex 2-orbifold, and holomorphic maps. A complex 2-orbifold $\Sigma$ is a
  Riemann surface $S$ marked by a finite set of marked points $\left\{ (q_1,
  m_1) \ldots, (q_n, m_n) \right\} $, where the $m'_i s$ are integers
  $\geqslant 2$.

  Let $X$ be a complex manifold, $f : X \rightarrow \Sigma$ a map. Let $x \in
  X, q = f (x) .$ Let $m \in \mathbbm{N}^{\ast}$ be the multiplicity of $q$,
  so that there exists \ an holomorphic map $u : D (0, r) \subset \mathbbm{C}
  \rightarrow (\Sigma, q)$ which is a ramified cover of order $m$ of a
  neighborhood of $q$. Then,$f$ is called holomorphic at x, if there exists a
  neighborhood $U$ of $q$ \ and a lift $\tilde{f} : U \rightarrow D$,
  holomorphic at $x$ such that $f = u \circ \tilde{f}$.
\end{definition}

\begin{definition}
  Fundamental group. Let $\Sigma = (S ; \left\{ (q_1, m_1) \ldots, (q_n, m_n)
  \right\})$be a $2 -$orbifold. Let $q \in S \backslash \left\{ (q_1, m_1)
  \ldots, (q_n, m_n) \right\}$. The fundamental group -in the sense of
  orbifolds- of $\Sigma$ at the point $p$ is the quotient \
  $\pi_1^{\tmop{orb}} (\Sigma, p) = \pi_1 (S \backslash \left\{ q_1, \ldots
  q_n \right\}) / \ll \gamma_i^{m_i} \gg$, where $\gamma_i$ is the class of
  homotopy (well defined up to conjugacy) of a small circle turning once
  around $q_i$, and $\ll \gamma_i^{m_i} \gg$ is the normal subgroup generated
  by all the conjugates of $\gamma _i^{m_i}.$ \ 
\end{definition}

\begin{example}
  (This is the main example, see $\left[ \tmop{Th} \right]$ chap. 13) Let
  $\Gamma \subset \tmop{PSL} (2, \mathbbm{R})$ be a uniform (discrete
  co-compact) lattice. The quotient $S = D / \Gamma$ of the unit disk by the
  action of $\Gamma$ is a Riemann surface. If $p \in D,$ its stabilizer is a
  finite hence cyclic subgroup of $\tmop{PSL} (2, \mathbbm{R}) .$ Modulo the
  action of $\Gamma$ there are only a finite set of points $\left\{ q_1,
  \ldots .q_n \right\} $with non trivial stabilizers of order $m_i$. The
  quotient orbifold is $\Sigma = (S ; \left\{ (q_1, m_1) \ldots, (q_n, m_n)
  \right\})$. One proves that $\Gamma = \pi_1^{\tmop{orb}} (\Sigma)$. An
  orbifold is called hyperbolic if it is obtained in this way ; an orbifold is
  hyperbolic if and only if its Euler characteristic $\chi^{\tmop{orb}}
  (\Sigma) = \chi (S) - \Sigma_{1 \leqslant i \leqslant n} (1 -
  \frac{1}{m_i})$ is non positive.
\end{example}

The following definition is useful to understand the structure of K\"ahler
groups (see $\left[ \tmop{ABCKT} \right])$.

\begin{definition}
  A K\"ahler manifold $X$ fibers if there exists a pair $(\Sigma, F)$ where
  $\Sigma = (S ; \left\{ (q_1, m_1) \ldots, (q_n, m_n) \right\})$ is a
  hyperbolic \ 2-orbifold, and $F : X \rightarrow \Sigma$ an holomorphic map
  with connected fibers$.$ Two such maps $F : X \rightarrow \Sigma, F' : X'
  \rightarrow \Sigma'$ are equivalent if the fibers of F and $F'$ are the same
  and images in $\Sigma$ and $\Sigma'$ of singular fibers have same order. In
  this case there exists an holomorphic isomorphism from $S$ to $S'$ which
  maps singular points of $S$ to singular points of $S'$ preserving the
  multiplicity.
\end{definition}

Let $\pi : X \rightarrow S$ be an holomorphic map from a compact complex
surface to a curve. If $q \in S$ is a singular value of $\pi$ , the analytic
set $\pi^{- 1} (q)$ can be decomposed in a finite union of irreducible sets,
$(D_i) .$ Away from a set of complex dimension $n - 2$ in $D_i$, hence of
complex codimension 2 in $X$, the map $p$ can by written $\pi (z_1, \ldots
.z_n) = z_1^{d_i}$, where $d_i \tmop{is} \tmop{the} \tmop{multiplicity}
\tmop{of} D_i$. The multiplicity of the fiber $\pi^{- 1} (q)$ is by definition
$m = \tmop{pgcd} (d_i)$. Let $\Sigma$ be the orbifold whose underlying space
is $S$, singular points are singular values of $\pi$ with corresponding
multiplicity.

\begin{lemma}
  $\pi : X \rightarrow \Sigma$ is holomorphic.
\end{lemma}

By construction, locally in the neighborhood of a point of $\pi^{- 1} (q), \pi
(x) = f_1^{d_1} \ldots .f_k^{d_k} + \tmop{cte}$, with \ $m | \tmop{pgcd} d_i$
$\Box$

The following finiteness theorem is well-known in the smooth case, and
implicit in the litterature at several places ; we give below a short proof
based on the hyperbolic geometry of hyperbolic orbifolds.

\begin{theorem}
  Let $X$ be a compact complex manifold. There exists, up to equivalence, a \
  finite set of pair $(\Sigma_i, F_i)$ where $\Sigma_i $is a complex
  hyperbolic $2$-orbifold, $F_i : X \rightarrow \Sigma_i$ is holomorphic with
  connected fibers $\Box$. 
\end{theorem}

Let us give a proof of this (well known) fact based on the
Kobayashi-hyperbolicity of a hyperbolic 2 orbifold $:$there exist no
holomorphic map from $\mathbbm{C}$ to an hyperbolic 2-orbifold as there exists
no holomorphic map from $\mathbbm{C}$ to the unit disk. Thus, by the Bloch
principle, \ as $X$ is compact there exists a uniform bound on the
differential of an holomorphic map $F : X \rightarrow \Sigma$. Therefore the
set of pairs $(F, \Sigma)$ is {\tmem{compact}} (two such orbifold are
$\varepsilon$-close if they are close for the Gromov-Hausdorff topology, i.e.
there exists a map between them which is isometric up to an error of
$\varepsilon$). But this compact space has only isolated points : if $F_1 : X
\rightarrow \Sigma_1$ is given, and, and the (Gromov-Hausdorff) distance of
$F$ to $F_1$ is smaller than the diameter of $\Sigma_1$ (for instance
$\leqslant 1 / 2 \tmop{diam} (X)$ where $X$ is endowed the Kobayashi
pseudo-metric) all the fibers of $F_1$ are send by $F$ inside a disk (or an
annulus in the case of the Margulis constant)therefore to a constant by the
maximum principle ; in other words $F$ factorizes through $F_1$ and induces an
isomorphism between $\Sigma$ and $\Sigma_1 . \Box$

\begin{remark}
  This proof shows that \ the number of pairs $(F, \Sigma)$ for a given
  complex manifold $X$ can be bounded by the Kobayashi diameter of $X$.

\end{remark}

The following is well known (see $\left[ \tmop{Si} 1 \right]$ $\left[
\right.$CKO$\left. \right]$) .

\begin{theorem}
  Let \ $F : X \rightarrow S$ by an holomorphic map with$_{}$ connected fibers
  from the complex manifold $X$ to a complex curve $S$. Let $\Sigma$ be the
  orbifold whose singular points are the singular values of $p$ and
  multiplicity the multiplicity of the corresponding fiber. Let $Y = F^{- 1}
  (b)$ be the fiber of a non singular point of $S.$ Let $\pi_1' (Y)$the image
  in $\pi_1 (X)$ of $\pi_1 (Y)$. One has the exact sequence

  \ \ \ \ \ \ \ \ \ \ \ \ \ \ \ \ \ \ \ $1 \rightarrow \pi'_1 (Y)
  \rightarrow \pi_1 (X) \begin{array}{l}
    \rightarrow
  \end{array} \pi_1^{\tmop{orb}} (\Sigma) \rightarrow 1$

  in particular the kernel of $\pi_{\ast} : \pi_1 (X) \rightarrow
  \pi_1^{\tmop{orb}} (B)$ is finitely generated.$\Box$
\end{theorem}

\subsection{Valuations.}

The next result is a reformulation of a fibration theorem of Gromov-Shoen
$\left[ \tmop{GS} \right] $and Simpson $\left[ \tmop{Si} 3 \right] $in terms
of the exceptional set in the sense of Bieri Neumann Strebel ; see also
$\left[ \tmop{De} \right]$ for a more general study of the BNS invariant of a
K\"ahler group, where $\omega \in H^1 (\Gamma, \mathbbm{R})$ rather than $H^1
(\Gamma, \mathbbm{Z})$).

\begin{theorem}
  Let $\omega \in H^1 (\Gamma, \mathbbm{Z})$. Then $\omega$ is exceptional iff
  there exist a hyperbolic orbifold $\Sigma$, an holomorphic map $F : X
  \rightarrow \Sigma$ such that $\omega \in F^{\ast} H^1 (\Sigma, \mathbbm{Z})
  .$
\end{theorem}

Let $\eta$ be a closed holomorphic $(1, 0)$ form whose real part is the
harmonic representative of $\omega .$ Let $\tilde{X}$ the universal cover of
$X$, and $F : \tilde{X} \rightarrow \mathbbm{R}$ a primitive of $\tmop{Re}
\eta$. From the definition (Remark 5) of $\mathcal{E}^1$ we know that $F
\geqslant 0$ is not connected ; \ $\left[ \tmop{Si} 3 \right] $applies. One
can also apply the \ proof of corollary 9.2 of $\left[ \tmop{GS} \right]$ to
the foliation defined by the complex valued closed $(1, 0)$ from whose real
part is the harmonic representative of $\omega$ .

To prove the converse (which will not be used), one remarks that for every $w
\in H^1 (\Sigma, \mathbbm{Z}),$its pull back to $H^1 (\Sigma, \mathbbm{Z})$ is
exceptional, as $\pi_1^{\tmop{orb}} (\Sigma)$ is hyperbolic, and the kernel of
$\pi_1^{\tmop{orb}} (\Sigma) \rightarrow \mathbbm{Z}$ cannot be finitely
generated. $\Box$

\subsection{The Green-Lazarsfeld set of a K\"ahler group.}

Let $K$ be a field. Recall that a character $\chi : \Gamma \rightarrow
K^{\ast}$ is called integral in the sense of Bass $\left[ \right.$Ba$\left.
\right]$ if $\chi (\Gamma) \subset O$, the ring of algebraic integers of $K.$

\begin{proposition}
  Let $X$ be a K\"ahler manifold, $\chi \in E^1 (\Gamma, K^{\ast})$ be a
  character. If $\chi$ is not integral, $X$ fibers over a 2-orbifold $\Sigma$
  such that $\chi \in F^{\ast} E^1 (\pi_1^{\tmop{orb}} (\Sigma), K^{\ast}$). 
\end{proposition}

Proof. Let $v$ be some valuation such that $\omega = v \circ \chi \neq 0$. \
Let $\Gamma$ acts on $T_v .$ By prop.6 this action is exceptional. Applying
Thm. 18 we get a pair $F, \Sigma$ such that $\omega \in F^{\ast} H^1 (\Sigma,
\mathbbm{Z}) .$ From the exact sequence of Theorem 17, we see that $\pi'_1
(Y)$ is a finitely generated normal subgroup of $\Gamma$ made up with elliptic
elements. As $\pi'_1 (Y)$ is {\tmem{finitely generated}}, the subtree of $T_v$
made up with fixed points of $\pi'_1 (F)$ is not empty. As $\pi'_1 (Y)$ is
normal, it is invariant by the action of $\Gamma$. Therefore the boundary of
this tree contains at least $3$ distinct elements. Thus acting on $P^1 (K)$
$\pi'_1 (Y)$fixes three different points and is the identity : $\pi'_1 (Y)
\subset \ker \rho$, and $\rho$ descends to some character on
$\pi_1^{\tmop{orb}} (\Sigma)$. $\Box$

The following proposition is a reformulation of a result by Beauville $\left[
\right.$Be$\left. \right] (\tmop{Cor} 3.6)$, it will be used to study the
cohomology class of $v \circ \chi,$ for the archimedian valuation $v (z) = \ln
| z |$ an $\chi : \Gamma \rightarrow \mathbbm{C}^{\ast}$ a character.

\begin{proposition}
  Let $X$ be a K\"ahler manifold, $\chi \in E^1 (\Gamma, \mathbbm{C}^{\ast})$
  be \ character. If $| \chi | \neq 1$, there exist an holomorphic map $F : X
  \rightarrow \Sigma$ \ from $X$ to a 2-orbifold $\Sigma$ such that $\chi \in
  F^{\ast} E^1 (\pi_1^{\tmop{orb}} (\Sigma), K^{\ast}$).

\end{proposition}

Combining propositions 19 and 20, we get the description of the GL set of a
K\"ahler manifold in terms of its fibering over hyperbolic 2-orbifolds. It
generalizes results by M. Green R. Lazarsfeld $\left[ \tmop{GL} \right]$, A.
Beauville $\left[ \tmop{Be} \right]$, C. Simpson $\left[ \tmop{Si} 2 \right]$,
F. Campana $\left[ \tmop{Ca} \right]$, R. Pink D. Roessler $\left[ \tmop{PR}
\right]$, who studied the case where the field $K$ is the field of complex
numbers.

\begin{theorem}
  Let $\Gamma$ be the fundamental group of a K\"ahler manifold X, $(F_i,
  \Sigma_i)_{1 \leqslant i \leqslant n}$ the family of fibration of $X$ over
  hyperbolic $2$-orbifolds. Let $K$ be a field of characteristic $p$ (if $p =
  0$, $K =\mathbbm{C})$, $\bar{F}_p \subset K$ the algebraic closure of $F_p$
  in K. Then $E^1 (\Gamma, K)$ is made with a finite set of torsion characters
  (contained in $E^1 (\Gamma, \bar{F}_p$) if $p > 0$) and the union of
  $F_i^{\ast} \tmop{Hom} (\pi^{\tmop{orb}}_1 (\Sigma_i), K^{\ast}) .$
\end{theorem}

\begin{proof}
  We shall prove that a character $\chi$ which is not in the union$ \bigcup
  F_i^{\ast} \tmop{Hom} (\pi^{\tmop{orb}}_1 (\Sigma_i), K^{\ast})$ must be a
  torsion character of bounded order. Let us fix such a character $\chi$.

  From theorem 8, we know that there exists a finite number of fields
  $K_{\nu}$ and characters $\xi_{\nu}$ such that $H^1 (\Gamma, \xi_{\nu}) \neq
  0$, and \ for every $\chi \in E^1 (\Gamma, K)$ there exists an index $\nu$
  for which $\ker \xi_{\nu} \subset \ker \chi$. If $\xi_{\nu}$ is not
  integral, there exists a 2-orbifold $\Sigma$ and a holomorphic map $F : X
  \rightarrow \Sigma$ such that $\ker F_{\ast} \supset \ker \xi_v$ : therefore
  $\ker F_{\ast} \supset \ker \chi$ and $\chi \in F^{\ast} E^1
  (\pi_1^{\tmop{orb}} (\Sigma)) .$

  Thus, \ as $\chi \nin \bigcup F_i^{\ast} \tmop{Hom} (\pi^{\tmop{orb}}_1
  (\Sigma_i), K^{\ast})$ \ $\chi$ is integral.

  Let us first discuss the case of positive characteristic. If $\xi_{\nu}$ is
  integral, then $\xi_{\nu} (\Gamma)$ is made with roots of unity of
  $K_{\nu}$. But $K_{\nu}$ is of finite transcendence degree over $F_p$ so
  admits only a finite number of roots of unity of degree $d_{\nu}$ (see
  $\left[ \tmop{Ba} \right]$ for instance). Therefore, $\chi$ is a torsion
  character of order $d$ dividing $d_{\nu}$.

  Suppose now that $\tmop{char} K = 0$, and $\xi_{\nu}$ is integral. Thus
  $K_{\nu}$ is a number field, and $\xi_{\nu} (\Gamma)$ is contained in the
  ring $O_{\nu}$ of integers of $\xi_{\nu} .$ If $| \xi_{\nu} | \neq 1$, or if
  one of its conjugates $\sigma (\xi_{\nu})$ has $| \sigma (\xi_{\nu}) | \neq
  1$, as $H^1 (\Gamma, \xi_{\nu}) \neq 0$ we know (prop. 20) \ that there
  exists a 2-orbifold $\Sigma$ and a holomorphic map $F : X \rightarrow
  \Sigma$ such that $\ker F_{\ast} \supset \ker \xi_v$ ; the previous argument
  apply and proves that $\chi \in F^{\ast} E^1 (\pi_1^{\tmop{orb}} (\Sigma))$.
  Therefore, $\chi$ must be a root of unity, by a theorem of Kronecker, of
  bounded degree $d$, as the degree of the $n$-th cyclotomic polynomial goes
  to infinity with $n$, and as $d$ divides the degree of $K_{\nu} .$ The rest
  of the argument is unchanged.

\end{proof}

Thus, the theorem 21 reduces the computation of $E^1 (\Gamma, K^{\ast})$ to
the case where $\Gamma$ is the fundamental group of a 2-orbifold.

\begin{proposition}
  Let $\Gamma = \pi_1^{\tmop{orb}} (\Sigma)$, for $\Sigma = (S ; (q_i, m_i)_{1
  \leqslant i \leqslant n})$ a hyperbolic \ 2-orbifold then,
  
  $E^1 (\pi^{\tmop{orb}}_1 (\Sigma_{}), K^{\ast}) = \tmop{Hom}
  (\pi_1^{\tmop{orb}} (\Sigma, K^{\ast}))$ unless $g = 1$ and for all $i$,
  $m_i \nequiv 0 (\tmop{char} K)$.
  
  \ \ \ \ \ \ \ \ \ \ \ \ If $g = 1$ and for all $i$ $m_i \nequiv 0
  (\tmop{char} K)$, \ $E^1 (\pi^{\tmop{orb}}_1 (\Sigma_{}), K^{\ast})$ is
  finite, made of torsion characters.
  
  \ \ \ \ \ \ \ \ 
\end{proposition}

Let $\chi : \pi^{\tmop{orb}} (\Sigma) \rightarrow K^{\ast}$ be a
representation. If $\chi = 1, H^1 (\pi_1^{\tmop{orb}} (\Sigma), K^{\ast}) =
\tmop{Hom} (\pi_1^{\tmop{orb}} (\Sigma), K^{\ast}) \neq 0$. If $g >
1,$consider a simple closed curves on $S$ such that $c_{}$ are homologous to
0, which separated $S$ in two compact surface of positive genus $S_1, S_2$,
with common boundary $c$ and such that all singular points are in $S_2$; if $g
= 1$ consider a curve \ $c$, which bounds a disk $\bar{D}$ on $S$ containing
all singular points $q_i$, and let $S_1 = S \backslash \tmop{int} (D)$ be the
other component. One consider a representation $\chi : \pi_1^{\tmop{orb}}
(\Sigma) \rightarrow K^{\ast}$, and note that $\chi (c) = 1$ as $c$ is
homologous to $0$. We think of $\chi$ as a local system on $\Sigma$ and we
will use a Mayer Vietoris exact sequence.

First note that if $\chi_{} |_{\pi^{} (S_1)}$ and $\chi |_{\pi^{\tmop{orb}}
(\Sigma_2)}$ are not 1, then $H^1 (\pi_1^{\tmop{orb}} (\Sigma), K^{\ast}) \neq
0$ : let $x_0 \in K$, there exists a unique 1-cocyle $c$ such that $c (g) =
x_0 (1 - \chi (g))$ is $g \in \pi_1 (S_1)$ , $c (g) = 0$ if $g \in S_2$.

If $\chi |_{S_1} = 1$, as $H^1 (S_{1,}, \partial S_1, K) = K^2,$one can find
a 1-cocyle $c$ whose restriction on $S_2$ or $D$ is 0, and restriction on $S^1
\tmop{is} \tmop{not} \tmop{trivial} .$

We are left to the case $\chi |_{S_2}$ or $\chi |_D = 1.$ If $g (S_2) > 0 H^1
(\pi_1^{\tmop{orb}} (\Sigma_2), C, K) \twoheadrightarrow K^{2 g}$ and the
previous argument apply. \

The remaining case is $g = 1,$ $\chi |_{\pi_1^{\tmop{orb}} (D)} = 1_{}$,
$\chi |_{\pi_1 (S_1)} \neq 1$. Note that in this case $\chi$ is a torsion
character. \ Furthermore, $H^1 ( \pi_1^{\tmop{orb}} (\Sigma_2), K) =$ $\left\{
(z_1, \ldots, z_n) \in K / m_i z_i = 0 \right\}$. This space is 0 unless $m_i
\equiv 0$ (char K) for some $i$. \ On the other hand, if $\rho |_{\pi_1 (S_1)}
\neq 0$ the homomorphism \ $H^1 ( \pi_1^{} (\Sigma_1), \rho) \rightarrow K$
which sends $\theta$ to $\theta (c)$ is an isomorphism. Using the exact
sequence of \ Mayer Vietoris, we see that \ $H^1 (\pi^{\tmop{orb}}_1
(\Sigma_{}), \chi) \neq 0$ if $g > 1$ or $g = 1$ and for some $i$, $m_i$
divides the characteristic of $K$. $\Box$

\section{Bibliography.}

$\left[ \tmop{ABCKT} \right]$ J. Amorous, M. Burger, K. Corlette, D.
Kotschick, D. Toledo , {\tmem{Fundamental groups of compact K\"ahler
manifolds}}. Mathematical Surveys and Monographs, 44. American Mathematical
Society, Providence, RI, 1996

$\left[ \tmop{AN} \right] $D. Arapura and M. Nori Solvable fundamental groups
of Algebraic Varieties and K\"ahler Manifolds. Compositio Math. 116, 173--188
(1999).

$\left[ \tmop{Ba} \right] $H. Bass. Groups of integral representation type
Pacific J. of Math. vol. 86 1, 1980.

$\left[ \tmop{Be} \right] $A. Beauville Annulation du $H^1$ pour les fibr\'es
en droites plats. Lecture Notes in Math. 1507, pp. 1--15, Springer-Verlag
(1992).

$\left[ \right.$Bi$\left. \right]$ Bieri, R. , Strebel R. Geometric invariant
for discrete groups. Preprint.

$\left[ \tmop{BG} \right]$ Bieri, Robert, Groves, J. R. J. The geometry of the
set of characters induced by valuations. J. Reine Angew. Math. 347 (1984),
168--195.

$\left[ \tmop{BNS} \right] $ R. Bieri, W. Neumann, R. Strebel, A geometric
invariant of discrete groups. Invent. Math. 90 (1987), no. 3, 451--477.

$\left[ \tmop{Bre} \right]$ E. Breuillard. On uniform exponential growth for
solvable groups. The Margulis Volume, Pure and Applied Math. Quarterly, to
appear.

$\left[ \tmop{Bro} \right] $Brown, Kenneth S. Trees, valuations, and the
Bieri-Neumann-Strebel invariant. Invent. Math. 90 (1987), no. 3, 479-504

$\left[ \tmop{Bru} \right]$ A. Brudnyi, \ Solvable quotients of K\"ahler
groups. Michigan Math. J. 51 (2003), no. 3, 477--490.

$\left[ \tmop{BT} \right] $Bruhat, F.; Tits, J. Groupes r\'eductifs sur un
corps local. \ Inst. Hautes Etudes Sci. Publ. Math. No. 41 (1972)

$\left[ \tmop{Ca} \right]$ Campana, F. Ensembles de Green-Lazarsfeld et
quotients r\'esolubles des groupes de K\"ahler. \ J. Algebraic Geom. 10
(2001), no. 4, 599--622. 32J27 (14F35)

$\left[ \tmop{CKO} \right]$ Catanese, F. Keum, J. Oguiso, K. Some remarks on
the universal cover of an open $K 3$ surface. Math. Ann. 325 (2003), no. 2,
279--286.

$\left[ \tmop{CS} \right]$ K. Corlette, C. Simpson, On the classification of
rank two representations of quasi-projective fundamental groups. Preprint,
Feb. 2007.

$\left[ \tmop{De} \right] $Delzant, T. L'invariant de Bieri, Neumann, Strebel
des groupes de K\"ahler, preprint 2006, soumis \`a Math. Annalen.

$\left[ \tmop{GL} \right]$ M. Green and R. Lazarsfeld Deformation theory,
generic vanishing theorems and some conjectures of Enriques, Catanese and
Beauville. Invent. Math. 90, 389--407 (1987).

$\left[ \tmop{GS} \right] $M. Gromov and R. Schoen, Harmonic maps into
singular spaces and p-adic superrigidity for lattices in groups of rank one,
Publ. Math. IHES 76 (1992), 165--246.

$\left[ \tmop{GR} \right] $Groves, J. R. J. Soluble groups with every proper
quotient polycyclic. Illinois J. Math. 22 (1978), no. 1, 90--95.

$ \left[ \tmop{Le} \right]$ G. Levitt $\mathbbm{R}$-trees and the
Bieri-Neumann-Strebel invariant. Publ. Mat. 38 (1994), no. 1, 195--202.

$\left[ \tmop{PR} \right] $R. Pink, D. Roessler, A conjecture of Beauville and
Catanese revisited. Math. Ann. 330 (2004), no. 2, 293--308.

$\left[ \tmop{Se} \right] $Serre, Jean-Pierre Arbres, amalgames, $\tmop{SL}
(2)$. \ Ast\'erisque, No. 46. Soci\'et\'e Math\'ematique de France, Paris,
1977. 189 pp.

$\left[ \tmop{Si} 1 \right]$ C. Simpson. The ubiquity of variations of Hodge
structure. Complex geometry and Lie theory (Sundance, UT, 1989), 329--348,
Proc. Sympos. Pure Math., 53, Amer. Math. Soc., Providence, RI, 1991.

$\left[ \tmop{Si} 2 \right] $C. Simpson Subspaces of moduli spaces of rank one
local systems. Ann. Sci. Ecole Norm. Sup. 26, 361--401 (1993).

$\left[ \tmop{Si} 3 \right] $C. Simpson. Lefschetz theorems for the integral
leaves of a holomorphic one-form. Compositio Math. 87 (1993), no. 1, 99--113.

$\left[ \tmop{Th} \right]$ Thurston, W.P. The geometry and topology of three
manifolds. Princeton 1978.

\end{document}